\documentclass[12pt]{amsart}
\usepackage{amssymb,amsmath,amsthm,amscd}
\usepackage[all]{xy}

\ifpdf\usepackage[colorlinks,pdftex, plainpages=false,colorlinks=false, hidelinks]{hyperref}
\else \usepackage[hypertex]{hyperref}
\fi

\hypersetup{pdftoolbar=false,plainpages=false,hyperindex=true,pdfdisplaydoctitle=true}

\numberwithin{equation}{section}

\numberwithin{equation}{section}
\newtheoremstyle{fancy1}{10pt}{10pt}{\itshape}{12pt}{\textsc\bgroup}{.\egroup}{8pt}{
}
\newtheoremstyle{fancy2}{10pt}{10pt}{}{12pt}{\itshape}{.}{8pt}{ }

\theoremstyle{fancy1}

\newtheorem*{prop*}{Proposition}

\newtheorem*{prob*}{Problem}

\newtheorem*{main*}{Theorem}

\newtheorem*{cor*}{Corollary}

\newtheorem*{problem*}{Problem}

\theoremstyle{fancy2}

\newtheorem*{rem*}{Remark}

\newcommand{\cref}[1]{Corollary~\ref{#1}}

\newcommand{\eps}{\varepsilon}

\newcommand{\RP}{\mathbb{R\mkern1mu P}}
\newcommand{\CP}{\mathbb{C\mkern1mu P}}
\newcommand{\HP}{\mathbb{H\mkern1mu P}}
\newcommand{\CaP}{\mathrm{Ca}\mathbb{\mkern1mu P}^2}

\newcommand{\C}{{\mathbb{C}} }
\newcommand{\R}{{\mathbb{R}} }

\newcommand{\QH}{{\mathbb{H}} }
\newcommand{\Ra}{{\mathbb{Q}} }
\newcommand{\Sph}{\mathbb{S}}

\newcommand{\F}{\ensuremath{\operatorname{F}} }
\newcommand{\G}{\ensuremath{\operatorname{G}} }
\newcommand{\D}{\ensuremath{\operatorname{\D}} }
\newcommand{\SO}{\ensuremath{\operatorname{SO}} }

\newcommand{\Sp}{\ensuremath{\operatorname{Sp}} }
\newcommand{\U}{\ensuremath{\operatorname{U}} }
\newcommand{\SU}{\ensuremath{\operatorname{SU}} }
\newcommand{\Spin}{\ensuremath{\operatorname{Spin}} }

\newcommand{\T}{\ensuremath{\operatorname{T}} }
\renewcommand{\S}{\ensuremath{\operatorname{S}} }

\newcommand{\fg}{{\mathfrak{g}} }

\newcommand{\fh}{{\mathfrak{h}} }
\newcommand{\fm}{{\mathfrak{m}} }

\newcommand{\fsu}{{\mathfrak{su}} }

\def\con#1=#2(#3){#1 \equiv #2 \bmod{#3}}

\newcommand{\ml}{\langle}
\newcommand{\mr}{\rangle}

\newcommand{\tr}{\ensuremath{\operatorname{tr}} }
\newcommand{\diag}{\ensuremath{\operatorname{diag}} }

\newcommand{\rank}{\ensuremath{\operatorname{rk}} }

\newcommand{\Ad}{\ensuremath{\operatorname{Ad}} }

\DeclareMathOperator{\Id}{Id} 
 
\DeclareMathOperator{\spam}{span}

\newcommand{\no}{\noindent}

\renewcommand{\F}{\mathsf{F}}

\newcommand{\dif} {{\operatorname{d}}}

\begin{document}
\title{Geometrically formal homogeneous metrics of positive curvature}

\author{Manuel Amann}
\address{Karlsruher Institut f\"ur Technologie\\
76133 Karlsruhe, Germany}
\email{manuel.amann@kit.edu}
\author{Wolfgang Ziller}
\address{University of Pennsylvania: Philadelphia, PA 19104, USA}
\email{wziller@math.upenn.edu}
\thanks{The first author was supported by  IMPA and a research grant of the German Research Foundation DFG. The second author was supported by CAPES-Brazil, IMPA, the National Science
Foundation and the Max Planck Institute in Bonn.}

\begin{abstract}
A Riemannian manifold is called geometrically formal if the wedge product of harmonic forms is again harmonic, which implies in the compact case that the manifold is topologically formal in the sense of rational homotopy theory.
A manifold admitting a Riemannian metric of positive sectional curvature is conjectured to be  topologically formal.
Nonetheless, we show that among the homogeneous Riemannian metrics of positive sectional curvature a geometrically formal metric is either symmetric, or a metric on a rational homology sphere.
\end{abstract}

\bigskip

\maketitle

Compact manifolds of positive sectional curvature form an intriguing field of study. On the one hand, there are few known examples, and on the other hand the two main conjectures in the subject, the two Hopf conjectures, are still wide open.

\smallskip

 The most basic examples of positive curvature are the rank one symmetric spaces $\Sph^n$, $\CP^n$, $\HP^n$ and $\CaP$. Homogeneous spaces of positive curvature have been classified \cite{Be,BB}: there are the homogeneous flag manifolds due to Wallach, $W^6=\SU(3)/T^2$, $W^{12}=\Sp(3)/\Sp(1)^3$ and $W^{24}=\F_4/\Spin(8)$, the Berger spaces $B^7=\SO(5)/\SO(3)$ and $B^{13}=\SU(5)/\Sp(2)\cdot \S^1$, and the Aloff--Wallach spaces $W^7_{p,q}=\SU(3)/\diag(z^p,z^q, \bar z^{p+q})$ with \linebreak[4]$\gcd(p,q)=1$, $p\geq q> 0$. See e.g. \cite{Zi2} for a detailed discussion. Furthermore, we have  the biquotient examples due to Eschenburg \cite{E1,E2} and Bazaikin \cite{Baz} and the more recent cohomogeneity one example in \cite{De,GVZ}.

 \smallskip

All the known examples have the following remarkable properties: They are rationally elliptic spaces, i.e.~their rational homotopy groups $\pi_i(M)\otimes \Ra$ vanish from a certain degree $i$ on, and the even dimensional ones  have positive Euler characteristic.
For general simply-connected positively (or more generally non-negatively) curved manifolds, the Bott-Grove-Halperin conjecture claims  rational ellipticity, whilst  the Hopf conjecture asserts that their Euler characteristic is positive in even dimensions.

\smallskip

A (simply-connected) topological space is called (topologically) \emph{formal} if its rational homotopy type is a formal consequence of its rational cohomology algebra, or, equivalently in the case of a manifold, if its real cohomology algebra is weakly equivalent to its de Rham algebra. It is a classical result of rational homotopy theory that rationally elliptic spaces with positive Euler characteristic are formal, see e.g. \cite{FHT}. In fact, one easily sees that all known examples of positive curvature are formal, in even as well as in odd dimensions. It is thus natural to conjecture that positively curved manifolds are formal in general.

We mention  here that the situation is different in non-negative curvature. Homogeneous spaces $G/H$ naturally admit non-negative curvature and are rationally elliptic. If  $\rank H= \rank G$ they have positive Euler characteristic and are hence formal. On the other hand, in \cite{Ama12b, Kot11} one finds many examples of  non-formal homogeneous spaces.

\smallskip

Other classical examples of formal spaces are compact symmetric spaces and compact K\"ahler manifolds. In the case of symmetric spaces this simply follows from the fact that  harmonic forms are parallel. Thus in
\cite{Kot01} the notion of geometric formality was introduced: A Riemannian metric is
\emph{geometrically formal}
if wedge products of harmonic forms are again harmonic.  On a compact manifold the Hodge decomposition implies that a
manifold admitting a geometrically formal metric
 is also topologically formal. See \cite{Bae12} and \cite{Kot12} for some recent results on geometrically formal metrics in dimension 3 and 4, and \cite{Kot01,Kot03,Kot09,Kot11,OP,GN} for obstructions to geometric formality.

\smallskip

There are very few known examples of compact geometrically formal manifolds. In fact, to our knowledge they all belong to the following classes (see \cite{Kot01,Kot12,Kot11,Bae12})
\begin{itemize}
\item a Riemannian metric  all of whose harmonic forms are parallel,
\item a homogeneous metric on a manifold whose rational cohomology is isomorphic to the cohomology of $\Sph^p\times\Sph^q$ with either $p$ and $q$ both odd, or $p$ even and $q$ odd with $p>q$,
\item Riemannian products of the above and finite quotients by a group of isometries.
\end{itemize}

In the homogeneous case
geometric formality  is an obvious consequence of homogeneity, since harmonic forms must be invariant under the id component of the isometry group.
Homogeneous spaces which have the rational cohomology of  the product of  spheres are classified in \cite{Kr}, and in \cite{Kot11} it was shown that many of them are not homotopy equivalent to symmetric spaces. There are other metrics where all harmonic forms are parallel, besides the compact symmetric spaces. For example, any metric on a rational homology sphere or a K\"ahler metric on a rational $\CP^n$,  e.g. the twistor space of the quaternionic symmetric space $\G_2/\SO(4)$. If one allows the manifold not to be simply connected, there are many such examples, e.g. fake $\CP^2$ and $\CP^4$, see \cite{GN}, which are compact quotients of complex hyperbolic space. Although these spaces may be called topologically formal, this property usually has not the strong consequences known from rational homotopy theory unless the space is nilpotent. For quotients of products, as for example  $(M\times \R^n)/\Gamma$ with $M$ geometrically formal,  one simply observes that geometric formality is a local property.
\smallskip

It is the main result of this article that geometric formality is also rare in positive curvature:

\begin{main*}\label{theoA}
A homogeneous
geometrically
formal metric of positive curvature is either symmetric or a metric on a rational homology sphere.
\end{main*}

In \cite{Kot11},Theorem 25, it was shown that a metric on a non-trivial $\Sph^2$ bundle over $\CP^2$ cannot be formal. This includes the 6 dimensional flag manifold $W^6$, as well as the inhomogenous Eschenburg biquotient. We will show that any homogeneous metric on the other two flag manifolds $W^{12}$ and $W^{24}$ cannot be geometrically formal.
Of course, every  metric on a sphere is  geometrically formal, and every homogeneous metric on $\CP^{2n}$, $\HP^n$ and $\CaP$ is symmetric. The Berger space $B^7$ is geometrically formal as well, since it is a rational homology sphere. This leaves the Berger space $B^{13}$, the Aloff--Wallach spaces, and the homogeneous metrics on $\CP^{2n+1}$. For the Aloff--Wallach spaces, it was shown in \cite{Kot11} that the normal homogeneous metric is not geometrically formal, but this metric does not have positive curvature.

\smallskip

The recent example of positive curvature in \cite{De,GVZ} is a rational homology sphere and hence geometrically formal.
 It would be interesting to know if the only other known examples of positive curvature, i.e.  the 7 dimensional Eschenburg spaces and 13 dimensional Bazaikin spaces, can admit geometrically formal metrics. They have the same cohomology as $W_{p,q}$ and $B^{13}$, but our methods do not apply in this case since the isometry group is too small.

\smallskip

It would also be interesting to have some other examples of homogeneous spaces where some of the homogeneous metrics are geometrically formal. Although the methods in this paper can be used to check this, an example seems to be difficult to find. Any relationship in the cohomology ring puts strong restrictions on a geometrically formal metric.

\smallskip

To prove the theorem we use the elementary fact that the de Rham cohomology is isomorphic to the finite dimensional algebra of invariant forms, and hence closed and harmonic forms can be computed explicitly. The Berger space $B^{13}$ has the rational cohomology  of $\CP^2\times \Sph^9$ and the Aloff--Wallach space $W_{p,q}$ that of
$\Sph^2\times \Sph^5$. Hence there is a unique harmonic 2-form $\eta$ and to be geometrically formal implies that $\eta^3$ resp. $\eta^2$ must be 0 as a form. It turns out that even among the closed invariant forms there are none whose power is 0. For $W^{12}$ and $W^{24}$ there are relations in the cohomology ring that contradict geometric formality.  In the case of $\CP^{2n+1}$, the situation is more interesting. Here the condition is that $\eta^k$ must be harmonic for all $k$. But  already the harmonic 4-form changes with the metric and is the square of the harmonic 2-form only if the metric is symmetric. We point out that this metric is also almost K\"ahler, hence gives examples of such metrics which are not geometrically formal.

\vspace{3mm}

 In Section 1 we explain some background about homogeneous spaces and their cohomology. In Section 2 we deal with  $B^{13} $ and in Section 3 with the Aloff--Wallach spaces. In Section 4  we discuss $W^{12}$ and $W^{24}$, and in Section 5 $\CP^{2n+1}$.

\vspace{3mm}


\section{Preliminaries}\label{prelim}
We first discuss the methods we will use to prove our main theorem.

\no Let $M=G/H$ be a homogeneous space with $H$ the stabilizer group at a base point $p_0\in M$. Using a fixed biinvariant metric $Q$ on the Lie algebra $\fg$, we define an orthogonal splitting
$$
\fg=\fh+\fm  \ \ \text{with identification} \ \  \fm\simeq T_{p_0}M
$$
induced by the action fields $X^*$ via $X\in\fm\to X^*(p_0)$. The action of  $H$ on $T_{p_0}M$ then becomes the adjoint action ${\Ad_H}$ on $\fm$. Choose an $\Ad_H$ invariant and $Q$ orthogonal decomposition
$$
\fm=\fm_0\oplus  \fm_1\oplus\ldots \oplus\fm_k
$$
such that ${\Ad_H}_{|\fm_0}=\Id$ and ${\Ad_H}_{|\fm_i}$ is irreducible.
A metric of the form $$g={g_0}_{|\fm_0}+\lambda_1Q_{|\fm_1}
+\lambda_2Q_{|\fm_2}+\ldots+\lambda_kQ_{|\fm_k}$$ with $g_0$ an inner product on $\fm_0$ and $\lambda_i$  positive constants, is then a $G$ invariant metric on $M$. If the the $\Ad_H$ representations $\fm_i$, $i=1\ldots k$, are all inequivalent, every $G$ invariant metric has this form. If $\fm_i\simeq\fm_j$ the inner products between $\fm_i$ and $\fm_j$ can be described by $1,2$ or $4$ arbitrary constants, depending on wether the representations are real, complex, or quaternionic.
\smallskip

We will use the elementary fact that the DeRham cohomology is isomorphic to the cohomology of $G$ invariant forms. By homogeneity this in turn is isomorphic to
$$
H^*_{DR}(M)\simeq \left((\Lambda^* \fm)^H,d\right)
$$
of forms on $\fm$ invariant under the isotropy action. The differential of a $k$-form $\omega\in (\Lambda^k \fm)^H$ is again $H$ invariant and can be computed via the following formula:
\begin{equation}\label{dw}
 d\omega(u_1,\ldots,u_{k+1})=\sum_{i<j}(-1)^{i+j}
\omega([u_i,u_j]_{\fm},u_1,\ldots, \hat{u_i},\ldots, \hat{u}_j\ldots u_{k+1})
\end{equation}
for $u_i\in\fm$, where $[u_i,u_j]_{\fm}$ denotes the projection of $[u_i,u_j]$ into $\fm$. On $\Lambda^* \fm$ we use the inner product that makes $e_{i_1}\wedge e_{i_2}\wedge \ldots \wedge e_{i_r}$, $i_1< i_2 <\ldots<i_r$ into an orthonormal basis of $\Lambda^r \fm$  for any orthonormal basis $e_i$ of $\fm$. We denote the codiferential by $\delta$. Since $\ml \omega, d\eta\mr=\ml\delta \omega , \eta\mr$, a $G$ invariant form $\omega\in (\Lambda^r \fm)^H$ is harmonic if and only if
$$
d\omega=0 \ \ \text{ and}\ \ \ml \omega,d\eta\mr=0  \ \ \text{for all}\ \ \eta\in (\Lambda^{r-1} \fm)^H.
$$
This reduces the computation of the  DeRham cohomology and the harmonic forms to a finite dimensional purely Lie algebraic computation. The equations are in fact linear in the coefficients of $\omega$ in some basis, and quadratic in the coefficients of the metric.

In order to simplify the computation of the differentials $d\omega$ we observe the following. Let $e_i$ be a basis of $\fg$ where each basis vector lies either in $\fh$ or $\fg$ and denote, by abuse of notation, the dual basis of 1-forms again by $e_i$. Although the 1-forms $e_i$ are in general not $\Ad(H)$ invariant, and hence do not represent forms on $G/H$, we can nevertheless formally use \eqref{dw} to compute $de_i$. By using the product rule we can then compute $d\omega$ for any $r$-form $\omega=\sum e_{i_1}\wedge e_{i_2}\wedge\ldots\wedge e_{i_r}$, in particular the $\Ad(H)$ invariant forms.
To see that the formula in \eqref{dw} satisfies the product rule, observe that we could replace $[u_i,u_j]_{\fm}$ by $[u_i,u_j]$ since the $\fh$ component will evaluate to 0. But then it becomes the usual formula for the Lie algebra cohomology of $\fg$ and hence satisfies a product rule. Notice though that in this generality $d^2\omega$ does not have to be 0, unless $\omega$ is $H$ invariant. This is due to the fact that the proof that it vanishes, in the case of the Lie algebra cohomology, uses the Jacobi identity which does not hold if we take the $\fm$ component of all Lie brackets.

\bigskip

\section{The Berger space $B^{13}$}\label{sec01}

\bigskip

For the 13 dimensional Berger space $B^{13}=\SU(5)/\Sp(2)\cdot\S^1$, the embedding $\Sp(2)\cdot\S^1\subset\SU(5)$ is given by
 $\diag(zA,\bar{z}^4)$ for $A\in \Sp(2)\subset\SU(4)$ and $z\in\S^1$.
The manifold $B^{13}$ has the same DeRham cohomology as $\CP^2\times \Sph^9$. One can see this  for example by using the two homogeneous fibrations $$\S^1\to \SU(5)/\Sp(2)\to B^{13}  \ \text{and}\
\SU(4)/\Sp(2)\to \SU(5)/\Sp(2) \to \SU(5)/\SU(4)$$
 and the fact that $\SU(5)/\SU(4)=\Sph^9$, $\SU(4)/\Sp(2)=\SO(6)/\SO(5)=\Sph^5$ and that $B^{13}$ is simply connected. Thus there exists one harmonic 2-form $\eta$. Geometric formality requires $\eta^2$ to be harmonic, and $\eta^3=0$ on the level of forms. We actually do not need to explicitly compute the harmonic forms since we will show  that there are no closed invariant 2-forms $\omega$ with  $\omega^3=0$.

 \smallskip

To compute the invariant forms, we first make the following observations.
 $\Sp(n)$ acts on $\QH^n$ via matrix multiplication and $\Sp(n)\cdot \Sp(1)$ via $(A,q)v=Avq^{-1}$ for $A\in\Sp(n), q\in\Sp(1)$ and $v\in \QH^n$. It is well known that the algebra $\Lambda^*(\QH^n)^{\Sp(n)}$ of invariant forms is generated by the 3 symplectic forms, corresponding to the K\"ahler forms $\omega_I, \omega_J,\omega_k$ associated to the  3 complex structures coming from right multiplication  with $I,J,K\in \Sp(1)$. Right multiplication with $\Sp(1)$ acts on $\spam{\{I,J,K\} }\simeq\R^3$ via matrix multiplication by $\SO(3)$ under the two fold cover $\Sp(1)\to \SO(3)$. Thus if $\S^1\subset\Sp(1)$ is given by $e^{it}$, the algebra \begin{equation}\label{invariant}
 \Lambda^*(\QH^n)^{\Sp(n)\cdot\S^1} \text{ is  spanned by } \omega_I \text{ and its powers}.
 \end{equation}

 \bigskip

 From the inclusions $\Sp(2)\cdot\S^1\subset \SU(4)\cdot\S^1=\U(4)\subset\SU(5)$ it easily follows that the decomposition of $\fm$ into irreducibles under the action of $H=\Sp(2)\cdot \S^1$ is given by $\fm=V\oplus W$ with $\dim V=5$ and $\dim W=8$. On $V$ the factor $\S^1$ acts trivially and $\Sp(2)$ via matrix multiplication by $\SO(5)$ under the two fold cover $\Sp(2)\to \SO(5)$. On $W$ it acts via $(A,z)v=Avz^{-1}$ with $(A,z)\in \Sp(2)\times\S^1$. It follows that $\Lambda^*(V)^H$ is spanned by a 0-form and a 5-form, the volume form $v$, and $\Lambda^*(W)^H$ by $\omega_I$ and its powers by \eqref{invariant}. On the other hand
 $\Lambda^k(V)\otimes\Lambda^l(W)$ with $k,l>0$ contains no invariant forms since the $\S^1$ factor clearly acts non-trivially. Thus $(\Lambda\fm)^H$ is spanned by $v$ and $\omega_I$ as an algebra. Since there is only one invariant 2-form, $\omega_I$ must be harmonic, and similarly $\omega_I^2$ as well.  In order to obtain the cohomology ring of $B^{13}$, we need $dv\ne0$, but the only possibility, up to a multiple, is $dv=\omega_I^3$. This implies that $\omega_I^3\ne 0$ and hence no invariant metric can be geometrically formal.

\bigskip

\section{The Wallach spaces $W_{p,q}$}

\bigskip

Let $H=\S^1_{k}=\diag(e^{ik_1t},
e^{ik_2t},e^{ik_3t})\subset G=\SU(3)$ where $k_i$ are fixed integers with $\sum k_i=0$. The quotient $G/H=\SU(3)/\S^1_{k}$ was studied by Aloff--Wallach \cite{AW} who showed that it admits a homogeneous metric with positive sectional curvature if none of the $k_i$ is 0. We will show that in fact none of the homogeneous metrics, even in this special case, can be geometrically formal. This was shown to be the case for the metric induced by the biinvariant metric on $\SU(3)$, but this metric does not have positive curvature.

\smallskip
It is well known that the rational cohomology ring of  $\SU(3)/\S^1_{k}$ is that of $\Sph^2\times\Sph^5$, but they can be differentiated by a torsion group in $H^4$.
Thus there exists one harmonic 2-form $\eta$. To be geometrically formal, we need $\eta^2=0$  on the level of forms. As in the previous case, we will again show that there are no closed 2-forms with square 0, although the computation in this case is much more involved.
\smallskip

We choose the following basis for the Lie algebra of $\SU(3)$. To describe it, let $E_{ij}$ be the matrix which has a 1 in row i and column j, and 0 otherwise. Set
\begin{align*}
E_1&= E_{12}- E_{21}, & E_2&=-E_{13}+E_{31}, & E_3&=E_{23}-E_{32}\\
F_1&=iE_{12}+i E_{21}, & F_2&=i E_{13}+i E_{31}, & F_3&=i E_{23}+i E_{32}\\
H_1&=iE_{11}-i E_{22}, & H_2&=-i E_{11}+i E_{33}, & H_3&=i E_{22}-i E_{33}.
\end{align*}
 We also choose the biinvariant metric on $\fsu(3)$ given by $\ml A,B\mr =-\frac 12 \tr(AB)$ in which $E_i,F_i$ are orthonormal. Furthermore,  $H_i$ have unit length, are orthogonal to $E_i,F_i$,  and $\ml H_i,H_j\mr=-\frac 12$.
For the Lie brackets we  have:
\begin{align*}
[H_i,E_i]&=2 F_i, & [H_i,F_i]&=-2 E_i, &  [H_i,E_j]&=-F_j , & [H_i,F_j]&=E_j\\
[E_i,E_j]&=E_k, & [F_i,F_j]&=-E_k, &  [E_i,F_j]&=-F_k , & [E_i,F_i]&=2 H_i  \\
\end{align*}
where $i,j,k$ is a cyclic permutation of $1,2,3$. For the decomposition of $\fg$ we choose
$$\fg=\fh+ V_0 + V_1 +V_2+V_3$$
where
$$
V_0= \spam(\varepsilon),\, V_i= \spam(E_i,F_i) \ \text{for}\ i=1,\ldots,3.
$$
Here $\varepsilon$ needs to be $Q$-orthogonal to $\fh$ and of unit length, i.e.
$$
\varepsilon=\sum r_i H_i \ \text{with}\ \ (r_1-r_2)k_1+(r_3-r_1)k_2+(r_2-r_3)k_3=0
\ \text{and}\ \sum r_i^2 -\sum r_ir_j=1.
$$
The subspaces $V_i$ are invariant under the isotropy action by $H$. On $V_0$ it acts trivially  and on $V_i$ as:
\begin{align*}
\Ad(\diag(e^{i\theta_1},e^{i\theta_2},e^{i\theta_3})E_1&=(\theta_1-\theta_2)F_1, &
\quad \Ad(\diag(e^{i\theta_1},e^{i\theta_2},e^{i\theta_3})F_1=-(\theta_1-\theta_2)E_1\\
\Ad(\diag(e^{i\theta_1},e^{i\theta_2},e^{i\theta_3})E_2&=(\theta_3-\theta_1)F_2, &
\quad \Ad(\diag(e^{i\theta_1},e^{i\theta_2},e^{i\theta_3})F_2=-(\theta_3-\theta_1)E_2\\
\Ad(\diag(e^{i\theta_1},e^{i\theta_2},e^{i\theta_3})E_3&=(\theta_2-\theta_3)F_3, &
\quad \Ad(\diag(e^{i\theta_1},e^{i\theta_2},e^{i\theta_3})F_3=-(\theta_2-\theta_3)E_3
\end{align*}
where $\theta_i=k_i\cdot t$.

For the differential forms we use the basis of 1-forms dual to the basis $\varepsilon,E_i,F_i$ and by abuse of notation use the same letters.
Using \eqref{dw} and the above Lie brackets one easily obtains the following exterior derivatives of 1-forms:

\renewcommand{\arraystretch}{1.4}
\stepcounter{equation}

\begin{table}[ht]
\begin{center}
\begin{tabular}{c|c}
$w$ 	& $\dif w$ \\
\hline
$E_1$ & $\hspace{10pt} E_2\wedge E_3-F_2 \wedge F_3+s_1F_1\wedge \varepsilon $ \\
$E_2$ & $\hspace{10pt} E_3\wedge E_1-F_3\wedge F_1+s_2F_2\wedge \varepsilon$ \\
$E_3$ & $\hspace{10pt} E_1\wedge E_2-F_1\wedge F_2+s_3F_3\wedge \varepsilon$ \\
$F_1$ & $-E_2\wedge F_3-F_2\wedge E_3-s_1E_1\wedge \varepsilon$ \\
$F_2$ & $-E_3\wedge F_1-F_3\wedge E_1-s_2E_2\wedge \varepsilon$ \\
$F_3$ & $-E_1\wedge F_2-F_1\wedge E_2-s_3E_3\wedge \varepsilon$\\
$\varepsilon$ & $s_1 E_1\wedge F_1+s_2 E_2\wedge F_2+s_3 E_3\wedge F_3$
\end{tabular}
\end{center}
     \vspace{0.3cm}
     \caption{Differentials of one-forms}\label{1forms}
\end{table}
where
$$
s_i=2r_i-r_j-r_k \ \text{with}\ i,j,k \ \text{distinct}.
$$
In this computation we need to use the fact that
$$
(H_1)_\fm=Q(H_1,\varepsilon)\eps=Q\big(H_1,\sum r_iH_i\big)\eps=(r_1-\tfrac12 (r_2+r_3))
$$
and hence $(2H_i)_\fm=s_i\eps $.

\smallskip

As explained above, these one-forms are  not all well defined on $G/H$ but are useful for computing the exterior derivative of 2-forms via the product formula for forms.

\smallskip

The discussion now depends on the values of the 3 integers $k_i$ and we differentiate between 3 cases.

\subsection{All three $k_i$ are distinct.}

Assume that $H=\S^1_{k}=\diag(e^{ik_1t},
e^{ik_2t},e^{ik_3t})$ with all $k_i$ distinct. Since  the differences $k_i-k_j$ are then also all distinct, the actions of $H$ on $V_i$ are all non-trivial and inequivalent.
Hence an invariant metric depends on 4 parameters. The only invariant 1-form is $\eps$, and the only invariant 2-forms are the volume forms of $V_i$, i.e. $\omega_i=E_i\wedge F_i$.
Without having to compute which forms $\omega=\sum a_i\omega_i$ are closed, it is clear that an invariant metric cannot be geometrically formal since $\omega^2=0$ implies that $a_i=0$ for all $i$.
\begin{rem*}
One easily sees that the form $\omega=\sum a_i\omega_i$ is closed iff $\sum a_i=0$ and harmonic if in addition $\sum a_is_it_i^2=0$, where $t_i$ is the length of $E_i$ and $F_i$.
\end{rem*}

\subsection{One of the $k_i$ vanishes.}

Here we can assume, since cyclic permutations of the $k_i$ and changing the sign of all 3 does not change the homogeneous space, that $(k_1,k_2,k_3)=(0,-1,1)$. These are in fact precisely those Wallach spaces which do not admit an invariant metric with positive curvature. Nevertheless we will now show that even here there are no geometrically formal metrics. The action of $\Ad_H$ on $V_i$ is a rotation of speed 1 on $V_1$ and $V_2$, and speed 2 on $V_3$. Thus the space of invariant metrics is 6-dimensional. $\eps$ is still the only invariant one form, but now we have 5 invariant 2-forms:
$$
\omega_i=E_i\wedge F_i,\ i=1,2, 3, \ \text{and}\ \omega_4=E_1\wedge E_2+F_1\wedge F_2,\
\omega_5=F_1\wedge E_2-E_1\wedge F_2.
$$
For $\varepsilon$ we choose
$$
\varepsilon=(H_1+2H_2)/\sqrt{3}=\diag(-i,-i,2i)/\sqrt{3}  \ \text{and hence}\ (s_1,s_2,s_3)=(0,3,-3).
$$
From Table \ref{1forms} we easily obtain the exterior derivatives of the invariant 2-forms:

\stepcounter{equation}
\begin{table}[ht]
\begin{center}
\begin{tabular}{c|c}
$w$ 	& $\dif w$ \\
\hline
$\omega_1$ & $E_1\wedge E_2 \wedge F_3-E_1\wedge E_3 \wedge F_2+E_2\wedge E_3\wedge F_1-F_1\wedge F_2\wedge F_3$ \\
$\omega_2$ & $E_1\wedge E_2\wedge F_3-E_1\wedge E_3\wedge F_2+E_2\wedge E_3\wedge F_1-F_1\wedge F_2\wedge F_3 $ \\
$\omega_3$ & $E_1\wedge E_2\wedge F_3-E_1\wedge E_3\wedge F_2+E_2\wedge E_3\wedge F_1-F_1\wedge F_2\wedge F_3$ \\
$\omega_4$ & $ 3\,\omega_5\wedge \varepsilon$ \\
$\omega_5$ & $ 3\,\omega_4\wedge \varepsilon$
\end{tabular}
\end{center}
     \vspace{0.3cm}
     \caption{Differentials of 2-forms for  $(k_1,k_2,k_3)=(0,-1,1)$ }
\end{table}

Thus the 2-form $\omega=\sum a_i\omega_i$ is closed if and only if $\sum a_i=0$ and $a_4=a_5=0$, and as in the previous case it follows that $\omega^2=0$ implies $a_i=0$ for all $i$.

\subsection{Two of the $k_i$ are equal.}

Up to permutations, we can assume that $(k_1,k_2,k_3)=(-2,1,1)$. Thus $\Ad_H$ acts with speed 3 on $V_1$ and $V_2$, but with opposite orientation, and trivially on $V_3$ and $V_0$. The metric is thus arbitrary on $V_0\oplus V_3$. Since the action on $V_1$ and $V_2$ are also equivalent, an invariant metric depends on 10 parameters. Now the  invariant 1-forms are $\varepsilon,\ E_3$ and $F_3$,
and the invariant 2-forms are:
$$
\omega_i=E_i\wedge F_i,\ i=1,2, 3, \  \omega_4=E_1\wedge E_2-F_1\wedge F_2,\
\omega_5=F_1\wedge E_2+E_1\wedge F_2,\ \omega_6=E_3\wedge \varepsilon,\ \omega_7=F_3\wedge \varepsilon.
$$
For $\varepsilon$ we choose
$$
\varepsilon=H_3  \ \text{and hence}\ (s_1,s_2,s_3)=(-1,-1,2).
$$

The differentials for the invariant 2-forms are:

\stepcounter{equation}
\begin{table}[ht]
\begin{center}
\begin{tabular}{c|c}
$w$ 	& $\dif w$ \\
\hline
$\omega_1$ & $E_1\wedge E_2 \wedge F_3-E_1\wedge E_3 \wedge F_2+E_2\wedge E_3\wedge F_1-F_1\wedge F_2\wedge F_3$ \\
$\omega_2$ & $E_1\wedge E_2\wedge F_3-E_1\wedge E_3\wedge F_2+E_2\wedge E_3\wedge F_1-F_1\wedge F_2\wedge F_3 $ \\
$\omega_3$ & $E_1\wedge E_2\wedge F_3-E_1\wedge E_3\wedge F_2+E_2\wedge E_3\wedge F_1-F_1\wedge F_2\wedge F_3$ \\
$\omega_4$ & $ -2( E_1 \wedge F_1 +E_2 \wedge F_2) \wedge F_3+2\omega_5 \wedge \varepsilon$ \\
$\omega_5$ & $-2( E_1 \wedge F_1 +E_2 \wedge F_2) \wedge E_3-2\omega_4 \wedge \varepsilon $ \\
$\omega_6$ & $( E_1 \wedge F_1 +E_2 \wedge F_2) \wedge E_3+\omega_4 \wedge \varepsilon$\\
$\omega_7$ & $( E_1 \wedge F_1 +E_2 \wedge F_2) \wedge F_3-\omega_5 \wedge \varepsilon$
\end{tabular}
\end{center}
     \vspace{0.3cm}
     \caption{Differentials of 2-forms for  $(k_1,k_2,k_3)=(-2,1,1)$ }
\end{table}
\no Thus a 2-form $\omega=\sum a_i\omega_i$ is closed if and only if
$$
a_1+a_2+a_3=0,\ -2a_4+a_7=0,\ -2a_5+a_6=0,
$$
 in other words, $a_1+a_2+a_3=0, a_7=2a_4$ and $a_6=2a_5$. This leaves us with a $4$-dimensional space of closed forms in degree $2$. One easily sees  that the square of such a closed form is 0 iff  all $a_i$ vanish.

This finishes the proof for the Aloff--Wallach spaces.

 \bigskip

\section{The flag manifolds $W^{12}$ and $W^{24}$}\label{flag}

\bigskip

The cohomology ring of the the 3 flag manifolds $W^{6}$, $W^{12}$ and $W^{24}$ is well know, and can be computed by using Borel's method for the cohomology ring of a homogeneous space $G/H$, see e.g. \cite{Bo1,Bo2}. In our case this is particularly simple since $\rank H=\rank G$ and since we can restrict ourselves to real coefficients.

The result is that it is generated by 3 elements $a_1,a_2,a_3\in H^{2k}(M,R)$, where $k=1,2,4$ for the 3 different flag manifolds. The relationships come from the Weyl group invariant polynomials, i.e., the symmetric polynomials in $a_i$ vanish. If we choose the generators $x=a_1+a_2$ and $y=a_1-a_2$ the cohomology ring is:
$$
H^*(M,R)=\{x,y\mid x^3=0, y^2=-3x^2\}
$$
with basis $x,y$ in dimension $2k$, as well as $x^2,xy$ in dimension $4k$, and the fundamental class $y^3$ in  dimension $6k$. The two relationships $ x^3=0$ and $ y^2=-3x^2$ put strong restrictions on a geometrically formal metric. For  $W^6$, the method in \cite{Kot11} used the fact that $y$ must be a symplectic form, whereas $x$ has a kernel, contradicting $ y^2=-3x^2$. This proof does not seem to work when $k>1$. Instead, we restrict ourselves to homogeneous metrics and use the algebra of invariant forms.

For all three flag manifolds $G/H$ we have the splitting
$$
\fm=V_1\oplus V_2\oplus V_3
$$
into $\Ad_H$ irreducibles, with $\dim V_i=2k$. Using representation theory, one easily sees that there are no invariant forms in degree $<2k$. In degree $2k$ we clearly have the $\Ad_H$ invariant volume forms $\omega_i$ of the modules $V_i$. Some differential must be nonzero since $b_{2k}=2$.  For $\SU(3)/\T^2$ and $\Sp(3)/\Sp(1)^3$ we also have inner automorphisms (e.g $\Ad(E_{12}-E_{21})$ for $\SU(3)/\T^2$ ) which interchange the 3 modules $V_i$.  For $F_4/\Spin(8)$ we have the triality automorphism of $\Spin(8)$. This outer automorphisms of $\Spin(8)$ also extends to inner automorphisms of $F_4$, see e.g. \cite{WZ}, Theorem 3.2,
 and takes $V_1$ to $V_2$, $V_2$ to $V_3$, and $V_3$ to $V_1$.  Thus there exist diffeomorphisms of $G/H$ which interchange the volume forms $\omega_i$, which implies that $d\omega_i\ne 0$ for all $i$. By rescaling $\omega_i$ if necessary we can assume that $\omega=\sum a_i\omega_i$ is closed iff
$a_1+a_2+a_3=0$.

From the description of the forms $\omega_i$ it is also clear that $\omega_i^2=0$, that $\omega_i\wedge \omega_j$, $i<j$ are linearly independent, and that $vol=\omega_1\wedge \omega_2\wedge \omega_3$ is a volume form. Thus $\omega^3=6a_1a_2a_3 vol$ is nonzero iff all three $a_i$ are nonzero. Hence $x$ must be one of 3 forms, depending which $a_i$ vanishes.  Assume that say $a_3=0$ and hence $x=\omega_1-\omega_2$ up to a multiple.
Then $y=\sum a_i\omega_i$ for some nonzero $a_i$ with $\sum a_i=0$. The $2k$ dimensional classes $x$ and $y$ are the only closed invariant forms and are hence harmonic.  But then the relation  $y^2=-3x^2$ in cohomology must also hold on the level of $4k$-forms. Since
$$
x^2=-2\omega_1\wedge \omega_2 \ \text{and}\ y^2=2a_1a_2\omega_1\wedge \omega_2
+2a_1a_3\omega_1\wedge \omega_3+2a_2a_3\omega_2\wedge \omega_3,
$$
 it follows that $a_1a_3=a_2a_3=0$, which implies that $y$ is a multiple of $x$. But this is not possible. This finishes the proof for the 3 Wallach flag manifolds.

 \bigskip

\section{The complex projective space $\CP^{2n+1}$}\label{CPn}

\bigskip

 For $\CP^{2n+1}=\SU(2n+2)/S(\U(2n+1)\U(1))$ it is well known that the set of homogeneous metrics can be described as follows, see e.g. \cite{Zi1}. First, observe that $\Sp(n+1)\subset \SU(2n+2)$ acts transitively on $\CP^{2n+1}$ with stabilizer $\Sp(n)\cdot\S^1$.
 From the inclusions $\Sp(n)\cdot\S^1\subset \Sp(n)\cdot\Sp(1)\subset\Sp(n+1)$ we obtain the twistor fibration $\Sph^2\to \CP^{2n+1}\to \HP^{n}$ and every homogeneous metric is a Riemannian submersion where one scales the metric induced by a biinvariant metric on $\Sp(n)$ with a factor $t$ on the fiber.

 \smallskip
Of course, on $\CP^{2n+1}$ we have only one harmonic 2-form $\alpha_H$, and a metric is geometrically formal if $\alpha_H^k$ are again harmonic for $k=1,2,\ldots 2n+1$. We will show that already for $\alpha_H^2$ this is only the case for the symmetric metric.
 \smallskip

 We  need to explicitly express the invariant form in some basis.
 We choose the embedding of $\Sp(n)$ in $\Sp(n+1)$ as the upper block embedding, i.e. the stabilizer of the last basis vector in its action on $\QH^n$. We first describe the basis of its orthogonal complement.
Recall that $E_{i,j}$ is the matrix which has a 1 in row i and column j, and 0 otherwise. Set
$$
 e_1=iE_{n+1,n+1},\ \ e_2=jE_{n+1,n+1},\ \ e_3=kE_{n+1,n+1}, \ \ Y_\alpha=E_{\alpha,n+1}-E_{n+1,\alpha}$$
 $$
 Y_{\alpha,1}=i E_{\alpha,n+1}+iE_{n+1,\alpha}
,\ \ Y_{\alpha,2}=j E_{\alpha,n+1}+jE_{n+1,\alpha} ,\ \ Y_{\alpha,3}=k E_{\alpha,n+1}+kE_{n+1,\alpha}$$
where $\alpha$ goes from $1$ to $n$.

Let $H=\Sp(n)\cdot\S^1\subset\Sp(n+1)$ where $\S^1=e^{it}\subset \Sp(1)$. Then the orthogonal complement  of $\fh$ in $\fg$ splits as
$$
\fm = V\oplus W=\C\oplus\QH^n \ \text{with}\ V=\spam(e_2,e_3)
\ \text{and} \ W=\spam( Y_\alpha,Y_{\alpha,1},Y_{\alpha,2},Y_{\alpha,3}), \ \alpha=1,\ldots n
$$
and $\diag(A,e^{it})\in H$ acts on $\fm$ as $(z,v) \to (e^{2it}z,Ave^{-it})$.
$H$ acts irreducibly on $V$ and $W$ and hence the metric depends on 2 parameters. We denote by $\ml \ , \ \mr_t$ the metric on $\fm$ where $e_i$ have length $t$ and the basis vectors in $W$ have length 1. Extended to a homogeneous metric on $G/H$, the symmetric metric then corresponds to $t=1$.

  \smallskip

For the $\fm$ component of the Lie brackets of vectors in $\fm$ we have:
$$
[e_2,e_3]_\fm =0\ ,\  [e_i,Y_\alpha]_\fm=-Y_{\alpha,i}\ ,\   [e_i,Y_{\alpha,i}]_\fm=Y_{\alpha}\ ,\ [e_i,Y_{\alpha,j}]_\fm=Y_{\alpha,k}\ ,\ [Y_{\alpha},Y_{\alpha,i}]_\fm=-2e_i$$
$$
 [Y_{\alpha,i},Y_{\alpha,j}]_\fm=2e_k\ ,\  [Y_{\alpha},Y_{\beta}]_\fm= [Y_{\alpha},Y_{\beta,i}]_\fm=[Y_{\alpha,i},Y_{\beta,i}]_\fm
=[Y_{\alpha,i},Y_{\beta,j}]_\fm=0$$
where $i,j,k$ is a cyclic permutation of $1,2,3$ and $\alpha,\beta$ are distinct.

As in the previous case, we first compute the differentials of 1-forms:

\renewcommand{\arraystretch}{1.4}
\stepcounter{equation}

\begin{table}[ht]
\begin{center}
\begin{tabular}{c|c}
$w$ 	& $\dif w$ \\
\hline
$Y_\alpha$ & $\hspace{10pt} e_2\wedge Y_{\alpha,2}+ e_3\wedge Y_{\alpha,3} $ \\
$Y_{\alpha,1}$ & $\hspace{10pt} e_2\wedge Y_{\alpha,3}- e_3\wedge Y_{\alpha,2}$ \\
$Y_{\alpha,2}$ & $\hspace{0pt} -e_2\wedge Y_{\alpha}\ \ + e_3\wedge Y_{\alpha,1}$ \\
$Y_{\alpha,3}$ & $\hspace{-7pt} -e_2\wedge Y_{\alpha,1}- e_3\wedge Y_{\alpha}$ \\
$e_2$ & $\hspace{10pt}\sum_\alpha ( -2Y_{\alpha}\wedge Y_{\alpha,2}+ 2Y_{\alpha,3}\wedge Y_{\alpha,1})$ \\
$e_3$ & $\hspace{10pt} \sum_\alpha ( -2Y_{\alpha}\wedge Y_{\alpha,3}+ 2Y_{\alpha,1}\wedge Y_{\alpha,2})$.
\end{tabular}
\end{center}
     \vspace{0.3cm}
     \caption{Differentials of one-forms on $\CP^{2n+1}$}\label{1formsCPn}
\end{table}

\smallskip

We now determine the $H$ invariant forms. Clearly, in $\Lambda^*(V)$ we only have the volume element $v=e_2\wedge e_3\in\Lambda^2(V)$. Recall that
the algebra $\Lambda^*(\QH^n)^{\Sp(n)}$ of invariant forms is generated by the 3 symplectic forms, corresponding to the K\"ahler forms $\omega_I, \omega_J,\omega_k$ associated to the  3 complex structures coming from right multiplication  with $I,J,K\in \Sp(1)$ on $W=\QH^n$. Thus $\Lambda^*(W)^{\Sp(n)}$ is generated by:
{ \fontsize{11}{15} \selectfont
$$
\omega_I=\sum_\alpha ( Y_{\alpha}\wedge Y_{\alpha,1}- Y_{\alpha,2}\wedge Y_{\alpha,3})\ , \omega_J=\sum_\alpha ( Y_{\alpha}\wedge Y_{\alpha,2}- Y_{\alpha,3}\wedge Y_{\alpha,1}), \ \omega_K=\sum_\alpha ( Y_{\alpha}\wedge Y_{\alpha,3}- Y_{\alpha,1}\wedge Y_{\alpha,2})
$$
}
and hence
$$
\Lambda^*(\fm)^{\Sp(n)}  \ \text{is generated by }\  e_1,\ e_2,\ \omega_I
,\ \omega_J,\ \omega_K.
$$
In this algebra we can identify the $H$ invariant forms by determining the action of the circle in $H$ and diagonalizing it via complexification.
On $e_2, e_3$ the circle $e^{it}\in \S^1\subset H$ acts via a rotation $R(2t)$ since it is given by conjugation. On the two-plane spanned by $Y_{\alpha},\ Y_{\alpha,1}$ one easily checks that it acts via  a rotation $R(-t)$ and on the two-plane spanned by $Y_{\alpha,2},\ Y_{\alpha,3}$ as  a rotation $R(t)$.  Hence it acts trivially on $\omega_I$, and on the two-plane spanned by  $\omega_J, \omega_K$ it acts as $R(2t)$.
This action is diagonal in the basis $e_2+ie_3,e_2-ie_3, \omega_I, \omega_J+i\omega_K,
\omega_J-i\omega_K$ and acts via $\theta^2 + (\theta^*)^2 + \Id+ \theta^2 + (\theta^*)^2$. Thus we obtain invariant forms, besides $\omega_I$, by taking real and imaginary parts of
$(e_2+ie_3)\wedge(e_2-ie_3)$ and $(  \omega_J+i\omega_K)\wedge
(\omega_J-i\omega_K)$ as well as $(e_2+ie_3)\wedge(\omega_J-i\omega_K)$.
This gives us the following basis for the invariant forms in low degrees:
$$
\Lambda^2(\fm)^H=\spam(v,\ \omega_I), \ \text{where }\  v=e_2\wedge e_3
$$
and
$$
\Lambda^3(\fm)^H=\spam(\beta_1,\beta_2)  \ \text{where }\ \beta_1=e_2\wedge \omega_J+e_3\wedge \omega_K \ \text{and }\
\beta_2=e_2\wedge \omega_K-e_3\wedge \omega_J
$$
and the invariant 4-forms
$$
\Lambda^4(\fm)^H=\spam(\omega_I^2,\ v\wedge \omega_I,\     \omega_J^2+  \omega_K^2).
$$

Notice that in the above language $de_2=- 2\omega_J $ and $de_3=- 2\omega_K $, and that we have the relations $v\wedge v=v\wedge\beta_1=v\wedge\beta_2=0$. Using Table \ref{1formsCPn} and the product formula,  one easily sees that:
$$
dv=2\beta_2,\ \ d\omega_I=2\beta_2,\ \ d\omega_J=2e_3\wedge\omega_I,\ \ d\omega_K=-2e_2\wedge\omega_I, \ \ d\beta_1=-2(\omega_J^2+  \omega_K^2 ) -4 v\wedge\omega_I.
$$
Thus $\alpha_H=v-\omega_I$ is the only closed 2-form, which is hence harmonic.

We now claim that $\alpha_H^2$ can only be harmonic for the symmetric metric.
For this, we compute the differentials of the 4-forms:
$$
d\omega_I^2=4\omega_I\wedge\beta_2,\ d(v\wedge \omega_I)=2\omega_I\wedge\beta_2,\    d(\omega_J^2+  \omega_K^2)=-4\omega_I\wedge\beta_2.
$$
Thus we have 2 closed 4-forms:
$$
\alpha_H^2=\omega_I^2-2v\wedge \omega_I,\  \ \text{and }\  \omega_J^2+  \omega_K^2+2 v\wedge \omega_I,
$$
and we need to determine which linear combination is harmonic. For this it needs to be orthogonal to the derivative of the invariant 3-forms, which is $d\beta_1$ since $d\beta_2=\frac12 d\omega_I^2=0$. Thus the 4-form is harmonic iff
$$
\ml a(\omega_I^2-2v\wedge \omega_I ) + b ( \omega_J^2+  \omega_K^2+2 v\wedge \omega_I),
(\omega_J^2+  \omega_K^2 ) +2 v\wedge\omega_I \mr_t=0.
$$
Observe that the inner products between the 3 symplectic forms are all the same, say equal to $L$, and that they are orthogonal to $v\wedge \omega_I$. Furthermore,  $\ml v\wedge\omega_I,v\wedge\omega_I\mr=
\ml v, v\mr\cdot \ml \omega_I,\omega_I\mr=
t^2L$. Thus we need
$$
2aL+4bL-4bt^2L=2L(a+2b(1-t^2))=0.
$$
But this implies that the only value of $t$ where $\alpha_H^2$ is harmonic is $t=1$, i.e. the symmetric metric.
\smallskip

This finishes the proof of our main Theorem.

\bigskip

We note that in the terminology from \cite{Na} we proved that a homogeneous metric on $\CP^{n}$ which is 2-formal, i.e., the product of harmonic 2-forms is again harmonic, is already symmetric.

We remark further that the metrics with positive sectional curvature are described as follows. For $B^{13}$ and $ \CP^{2n+1}$ we consider the fibrations $\Sph^2\to \CP^{2n+1}\to \HP^{n}$ and $\RP^5\to B^{13}\to
\CP^4$ and scale the fibers with $t$. The metric then has positive curvature iff $0<t<\frac43$.
For  the more complicated  description of the homogeneous positively curved metrics on $W_{p,q}$ see \cite{Pu}, for the ones on the flag manifolds see \cite{Va}, and for the ones on spheres \cite{VZ}.

\providecommand{\bysame}{\leavevmode\hbox
to3em{\hrulefill}\thinspace}


\begin{thebibliography}{10000}

\bibitem[Am]{Ama12b}
M.~Amann.
\newblock Non-formal homogeneous spaces.
\newblock {\em Math. Z.}, 274(3-4):1299--1325, 2013.

\bibitem[AW]{AW} S.~Aloff and N.~Wallach,
\emph{ An infinite family of 7--manifolds admitting positively
curved Riemannian structures}, Bull. Amer. Math. Soc. {\bf
81}(1975), 93--97.


\bibitem[Ba]{Bae12} C.~B\"ar,
\emph{ Geometrically formal 4-manifolds with nonnegative sectional curvature},
arXiv:1212.1325v2, 2012.

\bibitem[Baz]{Baz} Y.\ Bazaikin, \emph{On a family of $13$-dimensional
closed Riemannian manifolds of positive curvature}, Siberian Math.\
J., 37 (1996), 1068--1085.


\bibitem[BB]{BB} L.~B\'erard Bergery, {\em  Les vari\'et\'es
riemanniennes homog\`enes simplement connexes de dimension impaire
\`a courbure strictement positive}, J.\ Math.\ pure et appl.\ {\bf
55}\,(1976), 47--68.

\bibitem[Bo1]{Bo1} A.\ Borel,
\emph{Sur la cohomologie des espaces principaux et des espaces homogenes de groupes de Lie compacts}, Ann. of Math.\ \textbf{57} (1953), 115--207.

\bibitem[Bo2]{Bo2} A.\ Borel,
\emph{Sur l'homologie et la cohomologie des groupes de Lie compacts connexes}, Amer. J. of Math.\ \textbf{76} (1954), 273--342.

\bibitem[Be]{Be} M.\ Berger,
\emph{Les vari\'et\'es riemanniennes homog\`enes normales simplement
connexes \`a courbure strictement positive}, Ann.\ Scuola Norm.\
Sup.\ Pisa \textbf{15} (1961), 179--246.



\bibitem[De]{De}
O.~Dearricott, \emph{A 7-manifold with positive curvature}, Duke
Math. J. \textbf{158} (2011), 307--346..


\bibitem[E1]{E1}
J. H.~Eschenburg, \emph{New examples of  manifolds with strictly
positive curvature}, Invent. Math. \textbf{66} (1982), 469--480.

\bibitem[E2]{E2} J. H.~Eschenburg,
 {\em Freie isometrische Aktionen auf kompakten Lie-Gruppen
 mit positiv gekr\"ummten Orbitr\"aumen,}
 Schriftenr.~Math.~Inst.~Univ.~M\"unster {\bf 32} (1984).


\bibitem[FHT]{FHT}
Y.~Felix, S.~Halperin and J.-C..~Thomas, \emph{Rational homotopy theory}, Vol. 205 of \emph{Graduate Texts in Mathematics}, Springer-Verlag, New York, 2001.

\bibitem[GN]{GN}
J.-F. Grosjean and P.-A. Nagy,
\emph{On the cohomology algebra of some classes of geometrically
formal manifolds}, Proc. London Math. Soc. \textbf{98} (2009), 607--630.

\bibitem[GVZ]{GVZ}
K.~Grove, L.~Verdiani and W.~Ziller, \emph{An exotic $T_1\Sph^4$
with positive curvature}
Geom. Funct. Anal. \textbf{21} (2011), 499-524.

\bibitem[Ko1]{Kot01}
D. Kotschik,
\emph{On products of harmonic forms}, Duke Math. J. \textbf{107} (2001), 521--531.


\bibitem[Ko2]{Kot12}
D. Kotschik,
\emph{Geometric formality and non-negative scalar curvature}, arXiv:1212.3317, 2012.


\bibitem[KT1]{Kot03}
D. Kotschik and S. Terzic,
\emph{On formality of generalized symmetric spaces}, Math. Proc. Cambridge Phil. Soc. \textbf{134} (2003), 491--505.

\bibitem[KT2]{Kot09}
D. Kotschik and S. Terzic,
\emph{Chern numbers and the geometry of partial flag manifolds}, Comm. Math. Helv. \textbf{84} (2009), 587--616.

\bibitem[KT3]{Kot11}
D. Kotschik and S. Terzic,
\emph{Geometric formality of homogeneous spaces and biquotients}, Pacific J. Math. \textbf{249} (2011), 157--176.

\bibitem[Kr]{Kr}
L. Kramer, \emph{Homogeneous spaces, Tits buildings, and isoparametric hypersurfaces},
Mem. Amer. Math. Soc. \textbf{752}, American Mathematical Society, Providence, RI, 2002.

\bibitem[Na]{Na}
P.-A. Nagy, \emph{On length and product of harmonic forms
in K\"ahler geometry}, Math. Z. \textbf{254} (2006), 199-218.

\bibitem[OP]{OP}
 L.~Ornea and M.~Pilca, \emph{Remarks on the product of harmonic forms}, Pac. J. Math.  \textbf{250} (2011), 353–-363.

\bibitem[PW]{PW}
P.~Petersen and  F.~Wilhelm, \emph{An exotic sphere with positive
sectional curvature}, Preprint 2008.

\bibitem[Pr]{Pr}
G. Prasad and S.K. Yeung,
\emph{Arithmetic fake projective spaces and arithmetic fake Grassmannians},
Amer. J. Mathe. \textbf{131} (2009), 379--407.

\bibitem[P\"{u}]{Pu}
T.~P\"{u}ttmann, \emph{Optimal pinching constants of odd dimensional
homogeneous spaces}, Invent. math. \textbf{138}, (1999), 631–-684.

\bibitem[Va]{Va}
F.M. Valiev, \emph{Precise estimates for the sectional curvatures of homogeneous Riemannian
metrics on Wallach spaces}, Sib. Mat. Zhurn. \textbf{20} (1979), 248–-262.

\bibitem[VZ]{VZ} L.Verdiani - W.Ziller, \emph{Positively curved homogeneous metrics
on spheres}, Math. Zeitschrift, \textbf{261} (2009), 473–-488.


\bibitem[Wa]{Wa} N. Wallach, \emph{Compact homogeneous Riemannian
manifolds with strictly positive curvature}, Ann. of Math.,
\textbf{96} (1972), 277-295.

\bibitem[Zi1]{Zi1}
 W.~Ziller, \emph{Homogeneous Einstein metrics on Spheres and
projective spaces}, Math. Ann. \textbf{259} (1982), 351–-358.


\bibitem[WZ]{WZ} M.Wang-W. Ziller, \emph{On isotropy irreducible Riemannian
manifolds}, Acta. Math. 166 (1991), 223-261.

  \bibitem[Zi2]{Zi2} W.Ziller,
 \emph{Examples of Riemannian manifolds with nonnegative sectional curvature},
  in: Metric and Comparison Geometry,  Surv.  Diff. Geom. 11,
  ed. K.Grove and J.Cheeger, (2007), 63--102.



\end{thebibliography}
\end{document}